\documentclass[12pt, a4paper]{article}
\usepackage{amssymb}
\usepackage{latexsym}
\usepackage{amsbsy}

\newcommand{\bcen}{\begin{center}}     \newcommand{\ecen}{\end{center}}
\newcommand{\bay}{\begin{array}}      \newcommand{\eay}{\end{array}}
\newcommand{\beq}{\begin{eqnarray*}}      \newcommand{\eeq}{\end{eqnarray*}}

\def\char{\mbox{char}}
\def\gl{\mbox{gl.dim.}}

\def\rank{\mbox{rank}}

\def\Ext{\mbox{Ext}}

\def\dim{\mbox{dim}}
\def\mod{\mbox{mod}}

\def\Im{\mbox{Im}}
\def\Ker{\mbox{Ker}}

\def\Supp{\mbox{Supp}}

\def\la{\langle}
\def\ra{\rangle}
\begin{document}

\title{\large {\bf Hochschild cohomology of truncated quiver algebras
\footnote{Projects 10426014 and 10201004 supported by NSFC}}}

\author{\large Yunge Xu $^a$, Yang Han $^{b}$ and Wenfeng Jiang $^a$}

\date{\footnotesize a. Faculty of Mathematics $\&$ Computer Science,
Hubei University,\\ Wuhan 430062, P.R.China \hspace{5mm} E-mail:
xuy@hubu.edu.cn\\
b. Academy of Mathematics and Systems science, Chinese Academy of
Sciences,\\ Beijing 100080, P.R. China \hspace{5mm} E-mail:
hany@iss.ac.cn}

\maketitle

\begin{center}

\begin{minipage}{12cm}

{\bf Abstract}: For a truncated quiver algebra over a field of
arbitrary characteristic, its Hochschild cohomology is calculated.
Moreover, it is shown that its Hochschild cohomology algebra is
finite-dimensional if and only if its global dimension is finite
if and only if its quiver has no oriented cycles.

\medskip


\end{minipage}

\end{center}

\bigskip

\section*{Introduction}

\medskip

Let $k$ be a fixed field of arbitrary characteristic and $A$ a
finite-dimensional $k$-algebra. Let $A^e:=A^{op}\otimes_kA$ be the
{\it enveloping algebra} of $A$. The $n$-th {\it Hochschild
cohomology group} of $A$ is $H^n(A):=\Ext^n_{A^e}(A,A)$ and the
{\it Hochschild cohomology algebra} of $A$ is
$H^{\bullet}(A):=\coprod_{n \geq 0}H^n(A)$ where the
multiplication is induced by the Yoneda product for extensions
\cite{ML}.

\medskip

Let $Q=(Q_0,Q_1)$ be a finite connected quiver where $Q_0$ (resp.
$Q_1$) is the set of vertices (resp. arrows) in $Q$. Given a path
$p$ in $Q$, we denote by $o(p)$ and $t(p)$ the origin and the
terminus of $p$ respectively, and by $l(p)$ the length of $p$. For
$0 \leq j \leq l(p)$, denote by $^jp$ the subpath of $p$ with
length $j$ and origin $o(p)$. A path $a_1 \cdots a_n$ with $n \geq
1$ is called an {\it oriented cycle} of length $n$ if
$t(a_i)=o(a_{i+1})$ for $i=1,...,n-1$ and $t(a_n)=o(a_1)$. An
oriented cycle is called a {\it basic cycle} if it has no self
intersection.

\medskip

A {\it truncated quiver algebra} is an algebra $kQ/k^nQ$ where $n
\geq 2$ and $k^nQ$ denotes the ideal of $kQ$ generated by all
paths of length $\geq n$ (cf. \cite{C4}). It is further called a
{\it truncated basic cycle algebra} if $Q$ is a basic cycle. We
refer to \cite{ARS} for the theory of quivers and their
representations.

\medskip

Though Hochschild cohomology is theoretically computable for a
concrete algebra via derived functors, actual calculation for a
class of algebras is still very convenient, important and
difficult. So far the Hochschild cohomology was calculated for
hereditary algebras \cite{C2,Hap}, incidence algebras
\cite{C1,GS}, algebras with narrow quivers \cite{C2,Hap}, and so
on. For the Hochschild cohomology of truncated quiver algebras
partial results are known: Cibils calculated the Hochschild
cohomology of radical square zero algebras \cite{C3}. In case the
field $k$ is of characteristic zero, Zhang and Locateli calculated
the Hochschild cohomology of truncated basic cycle algebras and
truncated quiver algebras respectively \cite{Z0,Z,L}. For
truncated basic cycle algebras over a field $k$ of arbitrary
characteristic, Bardzell, Locateli and Marcos calculated their
Hochschild cohomology \cite{BLM}. In this paper, we complete the
calculation of the Hochschild cohomology of truncated quiver
algebras over a field of arbitrary characteristic. We shall employ
Locateli's strategy and do more and deeper analysis.

\medskip

In 1989, Happel asked the following question: Whether for a finite
dimensional algebra $A$ over an algebraically closed field $k$,
$\gl A < \infty$ if and only if the Hochschild cohomology algebra
of $A$ is finite-dimensional \cite{Hap}? Though Happel's question
has negative answer \cite{BGMS}, it does have some positive
answers: Avramov and Iyengar's result implies that Happel's
question has positive answer for commutative algebras \cite{AI}.
Cibils' result implies that Happel's question has positive answer
for radical square zero algebras \cite{C3}. Moreover, Locateli
proved that, in case the underlying field $k$ is of characteristic
0, the Hochschild cohomology algebra of a truncated quiver algebra
is finite-dimensional if and only if its quiver has no oriented
cycle. In this paper, we show by constructing explicitly some
Hochschild cocycles that for a truncated quiver algebra over a
field of arbitrary characteristic, its Hochschild cohomology
algebra is finite-dimensional if and only if its global dimension
is finite if and only if its quiver has no oriented cycles.

\medskip

\section{Hochschild cohomology groups}

\medskip

Let $A=kQ/k^nQ$, $n \geq 2$, be a truncated quiver algebra. If $S$
and $T$ are sets of paths in $Q$ we denote by $(S/\!\!/T)$ the set
$\{(p,q) \in S \times T | o(p)=o(q)$ and $t(p)=t(q)\}$. If $S$
(resp. $T$) is the set of paths of length $i$ (resp. $j$) then we
denote $(S/\!\!/T)$ by $(i/\!\!/j)$. Given a set $S$, $kS$ denotes
the $k$-vector space with basis $S$.

\medskip

If $Q$ is a basic cycle, Bardzell, Locateli and Marcos calculated
the Hochschild cohomology of $A$:

\medskip

{\bf Theorem 1.} (\cite[p. 1619]{BLM}) {\it Let $A=kQ/k^nQ$ with
$Q$ a basic cycle of length $e$ and $n=me+r$, $0 \leq r \leq e-1$.
Then we have} $$\dim_kH^0(A)=\left\{\begin{array}{ll} m,&\mbox{if $r=0$;}\\
m+e,&\mbox{if $r=1$;}\\ m+1,&\mbox{if $r \geq 2$.}
\end{array}\right.$$ $$\dim_kH^1(A)=\left\{\begin{array}{ll} m,&\mbox{if
$r=0,1$;}\\ m+1,&\mbox{if $r \geq 2$.}
\end{array}\right.$$ {\it and for} $i \geq 1$, $$\dim_kH^{2i}(A)=
\dim_kH^{2i+1}(A)=\left\{\begin{array}{ll} c_{n,e,i},&\mbox{if char$k \nmid n$ or $e \nmid (ni-n+1)$;}\\
c_{n,e,i}+1,&\mbox{if char$k|n$ and $e|(ni-n+1)$.}
\end{array}\right.$$ where $c_{n,e,i}:=|\{j|0 \leq j \leq n-2, j \equiv ri (\mod
e) \}|$.

\medskip

Thanks to Theorem 1, it remains to consider the case that $Q$ is
not a basic cycle. It follows from \cite[Lemma 2 and Lemma 6]{L}
that for all $i \geq 0$, the $i$-th Hochschild cohomology group of
$A$ is just the $i$-th cohomology group of the cochain complex
$$(M^{\bullet}, d^{\bullet}): \;\;\;\; 0 \rightarrow M^0
\stackrel{d^1}{\rightarrow} M^1 \stackrel{d^2}{\rightarrow} \cdots
\stackrel{d^{i-1}}{\rightarrow} M^{i-1}
\stackrel{d^i}{\rightarrow} M^i \stackrel{d^{i+1}}{\rightarrow}
M^{i+1} \stackrel{d^{i+2}}{\rightarrow} \cdots$$ where
$$\begin{array}{lcl} M^{2i}&=&k(0/\!\!/ni) \oplus \cdots \oplus
k(n-1/\!\!/ni),\\
M^{2i+1}&=&k(0/\!\!/ni+1) \oplus \cdots \oplus k(n-1/\!\!/ni+1),\\
d^{2i+1}&=&\left[\begin{array}{ccccc} 0&&&&\\ d^{2i+1}_0&0&&&\\
&d^{2i+1}_1&\ddots&&\\ &&\ddots&0&\\ &&&d^{2i+1}_{n-2}&0
\end{array}\right],\\
d^{2i}&=&{\footnotesize\left[\begin{array}{cccc} 0&0&\cdots&0\\
\cdots&\cdots&\cdots&\cdots\\ 0&0&\cdots&0\\ d^{2i}_0&0&\cdots&0
\end{array} \right]}\end{array}$$ with
$$d^{2i+1}_j: k(j/\!\!/ni) \rightarrow k(j+1/\!\!/ni+1), (p,q) \mapsto
\sum_{a \in Q_1}(ap,aq) - \sum_{a \in Q_1}(pa,qa)$$ for $0 \leq j
\leq n-2$ and $$d^{2i}_0: k(0/\!\!/n(i-1)+1) \rightarrow
k(n-1/\!\!/ni), (o(p),p) \mapsto \sum_{l(sq)=n-1}(sq,spq).$$
Therefore we have
$$\Ker \, d^{2i+1} = (\bigoplus_{j=0}^{n-2} \Ker \, d_j^{2i+1}) \oplus
k(n-1/\!\!/ni), \,\,\, \Im \, d^{2i+1} = \bigoplus_{j=0}^{n-2} \Im
\, d_j^{2i+1},$$
$$\Ker \, d^{2i} = \Ker \, d_0^{2i} \oplus \bigoplus_{j=1}^{n-1}
k(j/\!\!/n(i-1)+1), \,\,\, \Im \, d^{2i} = \Im \, d_0^{2i},$$
$$H^{2i+1}(A)=\Ker \, d^{2i+2}_0 \oplus\bigoplus^{n-1}_{j=1}(
k(j/\!\!/ni+1)/\Im \, d^{2i+1}_{j-1})$$ and
$$H^{2i}(A)=(k(n-1/\!\!/ni)/\Im \, d^{2i}_0) \oplus
\bigoplus^{n-2}_{j=0}\Ker \, d^{2i+1}_j.$$

\medskip

Now we recall some notations introduced by Cibils in \cite{C4}:
For $1 \leq j \leq n-2$ and $(p,q) \in (j/\!\!/ni)$. We say that
$(p,q)$ {\it starts together} (resp. {\it ends together}) if the
first (resp. last) arrow of $p$ and $q$ coincide. And we say that
$(p,q)$ {\it starts at a source} (resp. {\it ends at a sink}) if
$o(p)=o(q)$ is a source (resp. if $t(p)=t(q)$ is a sink).

Suppose that $(p,q)$ starts together, so $(p,q)=(ap',aq')$ for
some $a \in Q_1$, and moreover suppose that $(p,q)$ does not end
at a sink. A $+movement$ of $(p,q)$ is a couple $(p'b,q'b)$ where
$b \in t(p)Q_1$. If $(p,q)$ does not start together, or if $(p,q)$
end at a sink, then no $+$movement of $(p,q)$ is defined, and we
call $(p,q)$ a $+extreme$.

In analogous way, if $(p,q)$ ends together and $(p,q)$ does not
start at a source, then we define a $-movement$ of $(p,q)$. If
$(p,q)$ does not end together, or starts at a source, we say
$(p,q)$ is a $-extreme$.

Consider the equivalence relation on $(j/\!\!/ni)$ by declaring
$(p,q)$ {\it equivalent} to $(p',q')$ if there exists a finite
sequence of $+$ or $-$ movements transforming $(p,q)$ into
$(p',q')$. We call an equivalence class a $j$-{\it extreme} if all
its $+$extreme elements end at a sink and all its $-$extreme
elements start at a source.

\medskip

For notational convenience we denote $\{0$-extremes$\}=\emptyset$
and $\{(n-1)$-extremes$\}=(n-1/\!\!/ ni).$ One can check
Locateli's proofs step by step and obtain that the following
proposition and two lemmas hold for the underlying field $k$ of
arbitrary characteristic.

\medskip

{\bf Proposition 1.} (\cite[Proposition 3]{L}) {\it If $Q$ is not
a basic cycle, then} $$\dim_k \Ker \, d^{2i+1} =
\sum_{j=1}^{n-1}|\{j \mbox{-extremes}\}|$$ {\it and} $$\dim_k \Im
\, d^{2i+1} = \sum_{j=0}^{n-2}(|(j/\!\!/
ni)|-|\{j\mbox{-extremes}\}|).$$

\medskip

Let $x = \sum_{(p,q) \in (i/\!\!/j)}x_{(p,q)}(p,q) \in
k(i/\!\!/j)$. The {\it support} of $x$ is the set
$\Supp(x)=\{(p,q)|x_{(p,q)} \neq 0\}$.

\medskip

{\bf Lemma 1.} (\cite[Lemma 7]{L}) (Case $i=1$.) {\it Let
$(o(p),p) \in (0/\!\!/1)$ and} $(sq, spq) \in \Supp
(d^{2i}_0(o(p),p))$ {\it with $l(sq)=n-1$. Then $(sq, spq)$ is not
in the support of $d^{2i}_0(o(p'),p')$ for any other $p'\ne p$.}

\medskip

Let $p$ and $q$ be two paths in $Q$. We denote by $p \cap q$ the
set of the vertices that belong to both $p$ and $q$. Given a cycle
$p=a_1 \cdots a_j \cdots a_{n(i-1)+1}$ where $i \geq 2$ and $a_j
\in Q_1$ for $1 \leq j \leq n(i-1)+1$, the cycle $p \la a_1 \cdots
a_j \ra :=a_{j+1} \cdots a_{n(i-1)+1}a_1 \cdots a_j$ is called the
{\it rotate} of $p$ with origin $o(a_{j+1})$.

\medskip

{\bf Lemma 2.} (\cite[Lemma 8]{L}) (Case $i\geq 2$.) {\it Let
$(o(p),p) \in (0/\!\!/n(i-1)+1)$ with $i\geq 2$ and} $(sq, spq)
\in \Supp (d^{2i}_0(o(p),p))$ {\it with $l(sq)=n-1$. Let
$(o(p'),p') \in (0/\!\!/n(i-1)+1)$ and $(o(p),p) \ne (o(p'),p')$.
Then} $(sq, spq) \in \Supp (d^{2i}_0(o(p'),p'))$ {\it if and only
if $|sq \cap p| \geq 2$ and $p'$ is a rotate of $p$ with origin in
some vertex in $sq \cap p$.}

\medskip

{\bf Proposition 2.} {\it If $Q$ is not a basic cycle then the map
$d^{2i}_0: \; k(0/\!\!/n(i-1)+1)\rightarrow k(n-1/\!\!/ni)$ is
injective for all $i \geq 1$.}

\medskip

{\bf Proof.} Case $i=1$: Let $x=\sum_{(o(p),p) \in (0/\!\!/1)}
x_{(o(p),p)}(o(p),p) \in \Ker \, d^2_0$ where $x_{(o(p),p)} \in
k$. For any $(o(p),p) \in (0/\!\!/1)$ and $(sq, spq) \in \Supp
(d^2_0(o(p),p))$, by Lemma 1, we have that $(sq, spq)$ is not in
the support of $d^2_0(o(p'),p')$ for any $p' \neq p$. Since
$d^2_0(x)=0$, we must have $x_{(o(p),p)}=0$. So $x=0$ and $d^2_0$
is injective.

\medskip

Case $i \geq 2$: Define an equivalence relation on
$(0/\!\!/n(i-1)+1)$ by saying that $p$ is {\it equivalent} to $q$
if $p$ is a rotate of $q$. Denote by ${\cal C}$ the set of
equivalence classes. We can define a complete order in the set
$(0/\!\!/n(i-1)+1)$ as follows: First give a complete order $<'$
to the set of classes ${\cal C}$. Next, for each $c \in {\cal C}$,
choose $p=a_1 \cdots a_{n(i-1)+1} \in c$ and order the elements of
$c$ by the complete order $<''$ given by $p_1:=p<''p_2:=p \la a_1
\ra <''p_3:=p \la a_1a_2 \ra <''\cdots<''p_{n(i-1)+1}:=p \la a_1
\cdots a_{n(i-1)} \ra$. Then consider in $(0/\!\!/n(i-1)+1)$ the
lexicographic order: $(o(p),p)<(o(q),q) \, \mbox{ if }\,
\overline{p}<'\overline{q} \, \mbox{ or } \,
\overline{p}=\overline{q}\, \mbox{ and } \, p<''q,$ where
$\overline{p}$ and $\overline{q}$ denote the equivalence classes
of $p$ and $q$ respectively.

\medskip

We can also give a complete order to the set $(n-1/\!\!/ ni)$ as
follows: First, for any fixed equivalence class
$\overline{p}=\{p_1, p_2, \cdots,p_{n(i-1)+1}\} \in {\cal C}$, we
associate to it two subsets of $(n-1/\!\!/ni)$:
$${\cal B}_{\overline{p}}=\{(^{n-1}p_i,p_i \cdot \, ^{n-1}p_i) \in
(n-1/\!\!/ni)|p_i\in \overline{p}\}$$ and $${\cal
N}_{\overline{p}}=\{(sq,sp_iq) \in
(n-1/\!\!/ni)\setminus(\cup_{\overline{p} \in {\cal C}}{\cal
B}_{\overline{p}}) \mid p_i \in \overline{p}\}.$$ Next, consider
the complete order on ${\cal B}_{\overline{p}}$ induced by the
order in $\overline{p}$. Fix a complete order in ${\cal
N}_{\overline{p}}$ and assume that every element in ${\cal
B}_{\overline{p}}$ is smaller than that in ${\cal
N}_{\overline{p}}$. The order in ${\cal C}$ induces a complete
order in $(n-1/\!\!/ni)$ by assuming every element in ${\cal
B}_{\overline{p}} \cup {\cal N}_{\overline{p}}$ is smaller than
that in ${\cal B}_{\overline{q}} \cup {\cal N}_{\overline{q}}$ if
$\overline{p} < \overline{q}$ in ${\cal C}$. Finally, given any
complete order on $(n-1/\!\!/ ni)\backslash \bigcup_{c \in {\cal
C}}({\cal B}_c \cup {\cal N}_c)$ and require every element in
$\bigcup_{c \in {\cal C}}({\cal B}_c \cup {\cal N}_c)$ is smaller
than that in $(n-1/\!\!/ ni)\backslash \bigcup_{c \in {\cal
C}}({\cal B}_c \cup {\cal N}_c)$. Thus we defined a complete order
on $(n-1/\!\!/ ni)$.

\medskip

Under the ordered basis of $k(0/\!\!/n(i-1)+1)$ and
$k(n-1/\!\!/ni)$ respectively, the matrix of $d^{2i}_0$
can be written in the form $$\left[\begin{array}{cccc} B_1&&&\\
N_1&&&\\ &B_2&&\\ &N_2&&\\ &&\ddots&\\ &&&B_s\\ &&&N_s\\
0&0& \cdots &0
\end{array}\right]_{|(n-1/\!\!/ ni)| \times |(0/\!\!/n(i-1)+1)|}$$
where each submatrix $\left[\begin{array}{c}B_j\\
N_j \end{array}\right]$ corresponds to $d^{2i}_0|_{kc_j}$ and
${\cal C}=\{c_1, c_2, \cdots, c_s\}$.

\medskip

Given $(o(p),p) \in (0/\!\!/n(i-1)+1)$, Lemma 2 assures us that
$(^{n-1}p,p \cdot \, ^{n-1} p) \in \Supp (d^{2i}_0(o(p'),p'))$ if
and only if $p'$ is a rotate of $p$   and start at some vertex in
$(^{n-1}p) \cap p $. Let $e$ be the length of the smallest cycle
$w$ such that $p=w^h$ for some $h \geq 1$.

\medskip

(i). $n=2$ and $\char \, k \ne 2$: It follows from $e |
(n(i-1)+1)$ that $e \ne 2$. If $e=1$, then $p=a^{2i-1}$ for some
loop $a$ and $\overline{p}=\{p\}$. Since $Q$ is not a basic cycle,
$B_{\overline{p}}=(2)$ and $N_{\overline{p}}$ is in the form
$$\left[\begin{array}{c} 1\\ 1\\ \vdots \\ 1\end{array}\right].$$
Hence $$\rank
\left[\begin{array}{c} B_{\overline{p}}\\
N_{\overline{p}}\end{array}\right] = 1 = e.$$ If $e > 2$, the rank
of the matrix
$$B_{\overline{p}}=\left[\begin{array}{ccccc} 1 & 1& & & \\ &1&1 & &\\
& & 1 &\ddots & \\ & & &\ddots &  1 \\ 1& & & &
1\end{array}\right]_{e \times e}$$ is
$$\rank \, B_{\overline{p}}=\left\{\begin{array}{ll} e,&\mbox{if } \mbox{char }k \ne 2;\\
e-1,&\mbox{if }\mbox{char } k=2.
\end{array}\right.$$ Since $\char \, k \ne 2$, we have $\rank \,
B_{\overline{p}}=e$.

\medskip

Each block $B_{\overline{p}}, \, \overline{p} \in {\cal C},$ is of
full rank, so the map $d^{2i}_0$ is injective.

\medskip

(ii). $n>2$ and $\char \, k \nmid n$:  We write $n=me+r,$ with $0
\leq r <e$. It is not difficult to see that $(^{n-1}q, q \cdot \,
^{n-1}q)$ in ${\cal B}_{\overline{p}}$ appears in
$d^{2i}_0(o(p),p)$ with coefficient $m+1$ if $q=p_j$ with $1 \leq
j \leq r$, or $m$ if $q=p_j$ with $j>r$. Thus each matrix
$B_{\overline{p}}$ is of the form
$$B_{\overline{p}}=\left[\begin{array}{ccccc} m+1&m&m&\cdots&m+1\\
m+1&m+1&m&\cdots&m+1\\ m+1&m+1&m+1&\cdots&m+1\\
\vdots&\vdots&\vdots&&\vdots\\
m+1&m+1&m+1&\cdots&m+1\\
m+1&m+1&m+1&\cdots&m\\ m&m+1&m+1&\cdots&m\\ m&m&m+1&\cdots&m\\
\vdots&\vdots&\vdots&&\vdots\\
m&m&m&\cdots&m\\ m&m&m&\cdots&m+1\\ \end{array}\right]_{e \times
e}$$ which is a cyclic matrix, where in each column we have $r$
entries $(m+1)$ and $(e-r)$ entries $m$.

\medskip

Note that if $(o(p),p) \in (0/\!\!/n(i-1)+1)$ then one must have
$ni \equiv n-1$ (mod  $l(p)$). So $gcd(n,l(p))=1$ and thus
$gcd(e,r)=1$. By \cite[p. 1619]{BLM}, we have
$$\rank \, B_{\overline{p}}=\left\{\begin{array}{ll} e,&\mbox{if } \mbox{char } k \nmid n;\\
e-1,&\mbox{if }\mbox{char } k \mid n.
\end{array}\right.$$
Since char $k \nmid n$, we have $\rank \, B_{\overline{p}}=e.$

\medskip

Each block $B_{\overline{p}}, \, \overline{p} \in {\cal C},$ is of
full rank, so the map $d^{2i}_0$ is injective.

\medskip

(iii). $n > 2$ and char $k \mid n$, or $n=2$ and char $k=2$: Since
$Q$ is connected but not a basic cycle, without loss of
generality, we may assume that $Q$ contains an arrow $b$ which
does not appear in $w$ and $t(b)=o(p)$. Of course it is possible
that $o(b)$ belongs to the vertices of $w$. Moreover, we assume
$w=a_1a_2 \cdots a_e$ and $p=w^h$ for some $h \geq 1$ and
$n=me+r,$ with $0 \leq r < e.$ Clearly, $(b \cdot ^{n-2}p, b p
\cdot \, ^{n-2}p) \in {\cal N}_{\overline{p}}$ and $^{n-2}p=w^m
\cdot \, ^{r-2}p$.

\medskip

(1). If $m=0$ then $n-2=r-2$, thus $n-2 < e-2$. Hence
$$\begin{array}{lll} (b \cdot \, ^{n-2}p, b p \cdot \, ^{n-2}p)
&=&(b \cdot a_1 \cdots a_{n-2},
b(a_1 \cdots a_e)a_1 \cdots a_{n-2})\\
&=&(ba_1 \cdot a_2 \cdots a_{n-2},
ba_1(a_2 \cdots a_ea_1)a_2 \cdots a_{n-2})\\
&=& \cdots\\
&=& (ba_1 \cdots a_{n-2}, ba_1 \cdots a_{n-2}(a_{n-1} \cdots
a_ea_1 \cdots a_{n-2})).\end{array}$$ So $(b \cdot \, ^{n-2}p, bp
\cdot \, ^{n-2}p)$ appears only in $d^{2i}_0(o(p'),p')$ with
coefficient $1$  for $p'=p_j$ with $1 \leq j \leq n-1$, and with
coefficient $0$ otherwise. Thus the row corresponding to $(b\cdot
\, ^{n-2}p, bp \cdot \, ^{n-2}p)$ in the matrix $N_{\overline{p}}$
is of the form
$$(\underbrace{1,\,1,\,\cdots,\,1,}_{n-1} 0,\, 0,\, \cdots, 0)$$
whose row-sum is $n-1 \ne 0$. Therefore $$\rank
\left[\begin{array}{c} B_{\overline{p}}\\
N_{\overline{p}}\end{array}\right]=e.$$

\medskip

(2). If  $m \geq 1$ and $r \geq 2$ then $(b \cdot \, ^{n-2}p, bp
\cdot \, ^{n-2}p)$ appears in $d^{2i}_0(o(p'),p')$ for $p'=p_j$
with coefficient $m+1$ if $1 \leq j<r$, or $m$ if $r \leq j \leq
e$. So the row of the matrix $N_{\overline{p}}$ corresponding to
$(b \cdot \, ^{n-2}p, bp \cdot \, ^{n-2}p)$ is of the form
$$(\underbrace{m+1,\, m+1,\, \cdots, \, m+1}_{r-1},
\underbrace{m,\, \cdots, \, m}_{e-r+1})$$ whose row-sum is $n-1
\ne 0$. Thus $$\rank
\left[\begin{array}{c} B_{\overline{p}}\\
N_{\overline{p}} \end{array} \right]=e.$$

\medskip

(3). If $m \geq 1$ and $r=1$ then $^{n-2}p = w^{m-1} a_1a_2 \cdots
a_{e-1}$. So the coefficient of $(b \cdot \, ^{n-2}p, bp \cdot \,
^{n-2}p)$ in $d^{2i}_0(o(p'),p')$ for $p'=p_j$ is $m$ if $1 \leq j
\leq e$. Thus the row corresponding to $(b \cdot \, ^{n-2}p, bp
\cdot \, ^{n-2}p)$ in $N_{\overline{p}}$ is of the form
$$(\underbrace{m,\, m,\, \cdots, \, m}_{e})$$ whose row-sum is also
$(n-1)$ and hence nonzero. So $$\rank
\left[\begin{array}{c} B_{\overline{p}}\\
N_{\overline{p}}\end{array}\right] = e.$$

\medskip

(4). If $m \geq 1$ and $r=0$ then $^{n-2}p = w^{m-1}a_1a_2 \cdots
a_{e-2}$. The row corresponding to $(b \cdot \, ^{n-2}p, bp \cdot
\, ^{n-2}p)$ in $N_{\overline{p}}$ is of the form
$$( \underbrace{m,\, m,\,\cdots, \, m}_{e-1},m-1)$$ whose row-sum
is $n-1 \ne 0$. Therefore
$$\rank \left[\begin{array}{c} B_{\overline{p}}\\
N_{\overline{p}}\end{array}\right]=e.$$

\medskip

Anyway, each block $\left[\begin{array}{c} B_{\overline{p}}\\
N_{\overline{p}}\end{array}\right], \, \overline{p} \in {\cal C},$
is of full rank, so the map $d^{2i}_0$ is injective. \hfill
$\square$

\medskip

By Proposition 2, we immediately have:

\medskip

{\bf Proposition 3.} {\it If $Q$ is not a basic cycle then, for $i
\geq 1$, we have}
$$\dim_k \Ker \, d^{2i}=\sum_{j=1}^{n-1}|(j/\!\!/n(i-1)+1)|$$
{\it and} $$\dim_k \Im \, d^{2i}=|(0/\!\!/n(i-1)+1)|.$$

\medskip

Now we can complete the calculation of the Hochschild cohomology
of truncated quiver algebras over a field of arbitrary
characteristic as follows:

\medskip

{\bf Theorem 2.} {\it Let $A=kQ/k^nQ$ be a truncated algebra. If
$Q$ is not a basic cycle quiver then we have} $$H^0(A)=Z(A),$$
$$\dim_k H^1(A)=\dim_k
Z(A)-\sum_{j=0}^{n-1}|(j/\!\!/0)|+\sum_{j=1}^{n-1}|(j/\!\!/1)|$$
{\it and for} $i \geq 1$
$$\dim_k H^{2i+1}(A)=\sum_{j=1}^{n-1}|(j/\!\!/ni+1)|-
\sum_{j=0}^{n-2}(|(j/\!\!/ni)|-|\{j\mbox{-extremes}\}|)$$
$$\dim_k H^{2i}(A) = \sum_{j=1}^{n-1}|\{j\mbox{-extremes}\}|-|(0/\!\!/n(i-1)+1)|.$$
{\it where $Z(A)$ denotes the center of the algebra $A$.}

\medskip

{\bf Proof.} It is well-known that $H^0(A)=Z(A)$. Since $H^0(A) =
\Ker \, d^1$, we have $$\dim _k \Im \,d^1= \dim _k M^0 - \dim _k
\Ker \, d^1 = \sum^{n-1}_{j=0}|(j/\!\!/0)|- \dim _k Z(A).$$ By
Proposition 3, we have $$\dim _k H^1(A) = \dim _k \Ker \, d^2 -
\dim _k \Im \, d^1 = \sum^{n-1}_{j=1}|(j/\!\!/1)| -
\sum^{n-1}_{j=0}|(j/\!\!/0)| + \dim _k Z(A).$$ If $i \geq 1$ then
by Proposition 1 and Proposition 3 we have $$\begin{array}{lll}
\dim _k H^{2i+1}(A) &
 = & \dim _k \Ker \, d^{2(i+1)} - \dim _k \Im \, d^{2i+1}\\
& = & \sum_{j=1}^{n-1}|(j/\!\!/ni+1)| -
\sum_{j=0}^{n-2}(|(j/\!\!/ni)|-|\{j\mbox{-extremes}\}|)
\end{array}$$ and
$$ \begin{array}{lll}\dim _k H^{2i}(A) & = & \dim _k \Ker \, d^{2i+1} - \dim _k \Im
\, d^{2i}\\ & = &
\sum^{n-1}_{j=1}|\{j\mbox{-extremes}\}|-|(0/\!\!/n(i-1)+1)|.
\end{array}$$ \hfill{$\Box$}

\medskip

\section{Hochschild cohomolgy algebras.}

\medskip

The main result of this section is the following:

\medskip

{\bf Theorem 3.} {\it Let $A=kQ/k^nQ$ be a truncated quiver
algebra. Then the following conditions are equivalent:}

(1) {\it The algebra $H^*(A)$ is finite-dimensional;}

(2) $\gl A< \infty$;

(3) {\it $Q$ has no oriented cycles.}

\medskip

{\bf Proof.} (3) $\Rightarrow$ (2) $\Rightarrow$ (1): Clearly.

\medskip

(1) $\Rightarrow$ (3): If $Q$ has an oriented cycle then we take
the shortest oriented cycle $w=a_1a_2 \cdots a_e$ in $Q$. Thus the
vertices of $w$ must be different from each other.

Now we show that $(a_e,a_e(a_1 \cdots a_e)^{nr}) \in
k(1/\!\!/nre+1) \backslash \Im \, d^{2re+1}_0$, thus $d^{2re+1}_0$
is not surjective, hence $H^{2re+1}(A) \ne 0$, for all $r \geq 0$:

Assume on the contrary $(a_e,a_e(a_1 \cdots a_e)^{nr}) \in \Im \,
d^{2re+1}_0$. Then there is $\sum\limits_{i=1}^ex_i(v_i,(a_i\cdots
a_ea_1 \cdots a_{i-1})^{nr}) +
\sum\limits_{j=1}^sy_j(u_j,b_{j1}\cdots b_{j(nre)}) \in
k(0/\!\!/nre)$ with $\{b_{j1},...,b_{j(nre)}\} \nsubseteq
\{a_1,...,a_e\}$ such that $(a_e,a_e(a_1 \cdots a_e)^{nr})$

\noindent $\begin{array}{lll} &=&
d^{2re+1}_0(\sum\limits_{i=1}^ex_i(v_i,(a_i\cdots a_ea_1 \cdots
a_{i-1})^{nr}) +
\sum\limits_{j=1}^sy_j(u_j,b_{j1}\cdots b_{j(nre)}))\\
& =& \sum\limits_{i=1}^ex_i \sum\limits_{a \in
Q_1}((av_i,a(a_i\cdots a_ea_1 \cdots
a_{i-1})^{nr})-(v_ia,(a_i\cdots
a_ea_1 \cdots a_{i-1})^{nr}a))\\
&& + \sum\limits_{j=1}^sy_j\sum\limits_{a \in
Q_1}((au_j,ab_{j1}\cdots b_{j(nre)})-(u_ja,b_{j1}\cdots
b_{j(nre)}a))\\
& =& \sum\limits_{i=1}^ex_{i+1} (a_i,a_i(a_{i+1}\cdots a_ea_1
\cdots a_i)^{nr})\\
&&-\sum\limits_{i=1}^ex_i(a_i,(a_i\cdots
a_ea_1 \cdots a_{i-1})^{nr})a_i)\\
&&+\sum\limits_{i=1}^ex_i \sum\limits_{a \in
Q_1\backslash\{a_{i-1}\}}(av_i,a(a_i\cdots a_ea_1 \cdots
a_{i-1})^{nr})\\
&&-\sum\limits_{i=1}^ex_i\sum\limits_{a \in
Q_1\backslash\{a_i\}}(v_ia,(a_i \cdots a_ea_1 \cdots a_{i-1})^{nr}a)\\
&& + \sum\limits_{j=1}^sy_j\sum\limits_{a \in
Q_1}((au_j,ab_{j1}\cdots b_{j(nre)})-(u_ja,b_{j1}\cdots
b_{j(nre)}a)) \end{array}$ where $x_{e+1}:=x_1$.

\medskip

Let $W$ be the $k$-subspace of $k(1/\!\!/nre+1)$ generated by $e$
elements \linebreak $(a_i,a_i(a_{i+1} \cdots a_ea_1 \cdots
a_i)^{nr}), 1 \leq i \leq e.$ Then $k(1/\!\!/nre+1)=W \oplus
\overline{W}$ where $\overline{W}$ is the complement of $W$ in
$k(1/\!\!/nre+1)$. According to this decomposition we have
$(a_e,a_e(a_1 \cdots a_e)^{nr}) = \sum\limits_{i=1}^ex_{i+1}
(a_i,a_i(a_{i+1}\cdots a_ea_1 \cdots a_i)^{nr})
-\sum\limits_{i=1}^ex_i(a_i,(a_i\cdots a_ea_1 \cdots
a_{i-1})^{nr})a_i).$ Comparing the coefficients on two sides we
obtain

$$\begin{array}{lll} 1&=&x_1-x_e\\ 0&=&x_2-x_1\\
0&=&x_3-x_2\\ \cdots&\cdots&\cdots\\
0&=&x_e-x_{e-1}. \end{array}$$ It is a contradiction.
\hfill{$\Box$}


\medskip

\footnotesize

\begin{thebibliography}{99}
\bibitem{ML} Mac Lane S., Homology, Grundlehren 114, Third
corrected printing, Springer-Verlag, 1975.

\bibitem{C4} Cibils C., Rigidity of truncated quiver
algebras, Adv. Math., 1990, 79: 18-42.

\bibitem{ARS} Auslander M., Reiten I. and Smal\O \, S.O.,
Representation theory of artin algebras, Cambridge studies in
advanced mathematics 36, Cambridge university press, Cambridge,
1995.

\bibitem{C2} Cibils C., On the Hochschild cohomology of finite-dimensional
algebras, Comm. Algebra, 1988, 16: 645--649.

\bibitem{Hap} Happel D., Hochschild cohomology of
finite-dimensional algebras, Springer Lecture Notes in Math. 1404,
1989: 108--126.

\bibitem{C1} Cibils C., Cohomology of incidence algebras and
simplicial complexes, J. Pure Appl. Algebra, 1989, 56: 221--232.

\bibitem{GS} Gerstenhaber M. and Schack S.P., Simplicial homology
is Hochschild cohomology, J. Pure Appl. Algebra, 1983, 30:
143--156.

\bibitem{C3} Cibils C., Hochschild cohomology algebra of radical
square zero algebras, CMS Conf. Proc. 1998, 24: 93--101.

\bibitem{Z0}  Zhang P., Hochschild cohomology of truncated
algebras, Sci. China, Ser. A (in Chinese), 1994, 24: 1121-1125.

\bibitem{Z} Zhang P., Hochschild cohomology of truncated basic
cycle, Sci. China, Ser. A, 1997, 40: 1272-1278.

\bibitem{L} Locateli A.C., Hochschild cohomology of truncated
quiver algebras, Comm. Algebra, 1999, 27: 645-664.

\bibitem{BLM} Bardzell M.J., Locateli A.C. and Marcos E.N., On
the Hochschild cohomology of truncated cycle algebras, Comm.
Algebra, 2000, 28: 1615-1639.

\bibitem{BGMS}  Buchweitz R.O., Green E.L., Madsen D. and
Solberg \O., Finite Hochschild cohomology without finite global
dimension, Math. Res. Letters (in press).

\bibitem{AI} Avramov L.L. and Iyengar S., Gaps in Hochschild
cohomology imply smoothness for commutative algebras, Math. Res.
Letters (in press).

\end{thebibliography}
\end{document}